\documentstyle[12pt,a4]{article}
\begin{document}
\thispagestyle{empty}
\begin{center}
{\large\bf LIE ALGEBRA OF NONCOMMUTATIVE \\ INHOMOGENEOUS HOPF ALGEBRA}\\
\vskip2truecm M. Lagraa\footnote{e-mail : lagraa@elbahia.cerist.dz}and N. Touhami\footnote{e-mail : touhami@elbahia.cerist.dz}\\
\vskip1truecm Laboratoire de Physique Th\'eorique\\Universit\'e d'Oran Es-S\'enia, 31100, Alg\'erie\\
\vskip1truecm
\end{center}
\vskip3truecm 
\begin{abstract} We construct the vector space dual to the space of right-invariant differential forms construct from a first order differential calculus on inhomogeneous quantum group. We show that this vector space is equipped with a structure of a Hopf algebra which closes on a noncommutative Lie algebra satisfying a Jacobi identity.
\end{abstract}
\newpage
\section{Introduction}
It is well known that the underlying symmetry of the physical systems is governed by the Poincare group. In fact the most important observables in field theories such as the energy momentum and the angular momentum a
re related to the Lie algebra generators of the translations and the rotations. Then, it is especially interesting to study the noncommutative version of the Lie algebras of inhomogeneous Hopf algebras. In addition, due to the presence of the translations, one can hope to obtain new insight on the differential structure of the noncommutative space-time geeometry on which it is possible to construct consistently quantum field theories or a Einstein-Cartan gravity.\\Different approaches to construct Lie algebras of inhomogeneous quantum groups have been proposed [1-6]. Among them, there are which are based on the q-deformed Lie algebra of the Lorentz group augmented consistently by translation generators as q-deformed four-vector derivatives in Minkowski space [1-2]. Others approaches are based on the contraction of q-deformed anti-de sitter algebra ${\cal  U}(O(3,2))$ [4-5] or on the projective method where a bicovariant differential calculus on  inhomogeneous quantum groups as $IG
L_{q,r}(N)$ or $ISO_{q,r}(N)$ are obtained as quotients of those of $GL_{q,r}(N+1)$ or $ISO_{q,r}(N+1)$ with respect to a suitable Hopf ideal [6].\\
In this paper, one presents a quite different approach based on a bicovariant
differential calculus on an inhomogeneous Hopf algebra studied in our previous
paper [7]. This approach allows us to construct the vector space dual to the
right-invariant differential one-forms which is equipped with a Hopf algebra 
structure which closes on a quantum Lie algebra satisfying a quantum Jacobi 
identity. Throughout this paper we use the convention of those of Ref.[7-8].
\section{Differential Calculus On Inhomogeneous \\ 
Quantum Groups}
In this section, we start by recalling the definition and some basic properties of the 
inhomogeneous Hopf algebra ${\cal B}$ as well as the bicovariant first order 
differential calculus on it.\\ The quantum groups ${\cal G}$ are abstract 
objects described by corresponding Hopf algebra ${\cal B}$=$Poly({\cal G})$. 
For an inhomogeneous quantum group 
generated by an homogeneous quantum group ${\cal H}$ whose corresponding unital Hopf algebra is ${\cal A}=Poly({\cal H})$ and translations described by elements$p^{n},n=1,...,N$ corresponding to an irreducible representation 
$\Lambda^{n}_{~m}$ of ${\cal H}$ [8].\\We denote $\Delta : 
{\cal B}\rightarrow {\cal B} \otimes {\cal B}$, the comultiplication 
(algebra homomorphism), $\varepsilon : {\cal B} \rightarrow C$ the counit 
(character) and $S:{\cal B} \rightarrow {\cal B}$ the coinverse (algebra 
antihomomorphism) which satisfy the coassociativity $(\Delta \otimes id)\Delta =(id \otimes \Delta)\Delta$, $(\varepsilon \otimes id)\Delta = (id \otimes 
\varepsilon)\Delta = id$ and $m \circ (S \otimes id)\Delta = m \circ (id \otimes S)\Delta =I\varepsilon$, where $m : {\cal B} \otimes {\cal B} \rightarrow 
{\cal B}$ is the multiplication map ($m(a \otimes b)=ab$) and $I_{\cal B}=
I_{\cal A}=I$ is the unity of ${\cal B}$.\\These maps act on the generators of
${\cal B}$ as:
\begin{eqnarray*}
\Delta (\Lambda^{n}_{~m})=\Lambda^{n}_{~k} \otimes \Lambda^{k}_{~m} &,& \Delta(p^{n})=\Lambda^{n}_{~k} \otimes p^{k} + p^{n} \otimes I\\ \varepsilon(\Lambda^{n}_{~m}) = \delta^{n}_{m}  &,& \varepsilon (p^{n}) = 0, \\S(\Lambda^{n}_{~k})
\Lambda^{k}_{~m}=\Lambda^{n}_{~k}S(\Lambda^{k}_{~m})&=&\varepsilon ( \Lambda^{n}_{~m}),\\   S(\Lambda^{n}_{~k})p^{k}+S(p^{n})=\Lambda^{n}_{~k}S(p^{k})+p^{n}&=&\varepsilon (p^{n}).
\end{eqnarray*}
To construct a ${\cal B}$-bimodule of one-forms $\Gamma$, we consider an 
exterior bicovariant derivative $d:{\cal B} \rightarrow \Gamma$ 
($\Delta_{L}d=(id \otimes d)\Delta$ and $\Delta_{R}d=(d \otimes id)\Delta$) [7] from which we construct a right-invariant basis (Maurer-Cartan forms) of the bimodule $\Gamma$ over ${\cal B}$
\begin{eqnarray}
\Theta^{n}_{~m} = d \Lambda^{n}_{~k} S( \Lambda^{k}_{m})  &and~~\Pi^{n} = dp^{n} - \Theta^{n}_{~k}p^{k}
\end{eqnarray}
which transforms as
\begin{eqnarray}
\Delta_{R}( \Theta^{n}_{m})&=& \Theta^{n}_{m} \otimes I ,\\ \Delta_{R}
 ,\\\Delta_{L}( \Theta^{n}_{m})&=&\Lambda^{n}_{~\ell} S( \Lambda^{k}_{~m}) \otimes \Theta^{\ell}_{k} ,\\\Delta_{L}( \Pi^{n}) &=& \Lambda^{n}_{~k} \otimes \Pi^{k} + \Lambda^{n}_{~k} S(p^{\ell}) \otimes \Theta^{k}_{\ell}.
\end{eqnarray}
Folowing Ref.[7], we can state that there exist linear functionals $f^{n}_{~m}$, $f^{nk}_{~m}$, $f^{n}_{km}$ and $f^{n\ell}_{mk}$ $\in {\cal B}'$ (the dual vector space of ${\cal B}$) such that\\
\begin{eqnarray}
 \Pi^{n} a &=& (a \star f^{n}_{~k}) \Pi^{k} + (a \star f^{n\ell}_{~k}) \Theta^{k}_{~\ell} ,\\\Theta^{n}_{~m} a &=& (a \star f^{n}_{mk}) \Pi^{k} + (a \star f^{n\ell}_{mk}) \Theta^{k}_{~\ell}
\end{eqnarray}
and\\
\begin{eqnarray}
b \Pi^{n} &=& \Pi^{k}(b \star f^{n}_{~k} \circ S) + \Theta^{k}_{~\ell}(b \star f^{n\ell}_{~k} \circ S) ,\\b \Theta^{n}_{~m} &=& \Pi^{k} (b \star f^{n}_{mk} \circ S) + \Theta^{k}_{~\ell}(b \star f^{n\ell}_{mk} \circ S)
\end{eqnarray}
where the convolution product of a functional  $f \in {\cal B}'$ and an el
ement $a$ of ${\cal B}$ is defined as $(a \star f) = (f \otimes id)\Delta(a)$ or  $(f \star a) = (id \otimes f)\Delta(a)$ and the convolution product of two functionals is defined as $(f_{1}\star f_{2})(a) = (f_{1} \otimes f_{2})\Delta(a)$.\\From (6-7),we deduce
\begin{eqnarray}
f^{n}_{~m}(I)=\delta^{n}_{m},~~f^{nm}_{k\ell}(I)=\delta^{n}_{\ell}\delta^{m}_{k},~~f^{nm}_{~k}(I)=0&and~f^{n}_{km}(I)=0
\end{eqnarray}
and
\begin{eqnarray}
f^{n}_{~m}(ab)&=&f^{n}_{~k}(a)f^{k}_{~m}(b)+f^{nk}_{~\ell}(a)f^{\ell}_{km}(b),\\f^{nk}_{m\ell}(ab)&=&f^{np}_{mq}(a)f^{qk}_{p\ell}(b) + f^{n}_{mq}(a)f^{qk}_{\ell}(b),\\f^{n}_{m\ell}(ab)&=&f^{np}_{mq}(a)f^{q}_{p\ell}(b) + f^{n}_{mk}(a)f^{k}_{~\ell}(b),\\f^{nk}_{~\ell}(ab)&=&f^{n}_{~q}(a)f^{qk}_{~\ell}(b) + f^{np}_{~q}(a)f^{qk}_{p\ell}(b)
\end{eqnarray}
for any $a$ and $b \in {\cal B}$. From (10-14) we can see that the convolution product of different functionals generates an  algebra ${\cal B}^{0} \subset {\cal B}'$  which can be equipped with a Hopf algebra structure where th
\begin{eqnarray*}
\Delta'(f^{n}_{~m})&=& f^{n}_{~k} \otimes f^{k}_{~m} + f^{nk}_{~\ell} \otimes f^{\ell}_{km},\\\Delta'(f^{nk}_{m\ell})&=&f^{np}_{mq} \otimes f^{qk}_{p\ell} + f^{n}_{mq} \otimes f^{qk}_{~\ell},\\\Delta'(f^{n}_{m\ell})&=&f^{np}_{mq} \otimes f^{q}_{p\ell} + f^{n}_{mk} \otimes f^{k}_{~\ell},\\\Delta'(f^{nk}_{~\ell})&=&f^{n}_{~q} \otimes f^{qk}_{~\ell} + f^{np}_{~q} \otimes f^{qk}_{p\ell},
\end{eqnarray*}
 the counit is given by $\varepsilon'(f^{n}_{~m})=\delta^{n}_{m}\varepsilon$, $\varepsilon'(f^{nk}_{m\ell})=\delta^{n}_{\ell}\delta^{k}_{m}\varepsilon$ and $\varepsilon'(f^{n}_{mk})=\varepsilon'(f^{nk}_{~m})=0$ and the antipode is given by $S'(f)=f \circ S$. One can easily verify that these maps satisfy the Hopf algebra axioms
\begin{eqnarray*}
(\Delta' \otimes id)\Delta' &=& (id \otimes \Delta')\Delta',\\(\varepsilon' \otimes id)\Delta' &=& (id \otimes \varepsilon')\Delta' = id,\\m' \circ (S' \otimes id)\Delta' &=& m' \circ (id \otimes S')\Delta'= I_{{\cal B}^{0
}}\varepsilon'
\end{eqnarray*}
where $m'$ is the multiplication map defined as $m' \circ (f_{1} \otimes f_{2}) = f_{1} \star f_{2}$.\\In the following one assumes
\begin{eqnarray}
f^{nk}_{m\ell} = \tilde{f}^{k}_{~m} \star f^{n}_{\ell}  &~~f^{nk}_{~m} = \tilde{\eta}^{k} \star f^{n}_{~m}  &and~~f^{n}_{km} = \tilde{\eta}_{k} \star f^{n}_{~m},
\end{eqnarray}
where $\tilde{\eta}^{n}$, $\tilde{\eta}_{n}$ and $\tilde{f}^{n}_{~m} \in {\cal B}^{0}$. One deduces from (10) and (15) that 
\begin{eqnarray}
\tilde{\eta}^{n}(I) = 0 &,~~ \tilde{\eta}_{n}(I) = 0  &,and~~ \tilde{f}^{n}_{~m}(I) = \delta^{n}_{m}.
\end{eqnarray}
The fact that ${\cal A}$ is a Hopf sub-algebra of ${\cal B}$, we have [7]
\begin{eqnarray}
\tilde{\eta}_{n}(a) = 0 &,~~ a \in {\cal A}. 
\end{eqnarray}
From (6-7) and (17), one  can show [7] that the commutation rules between the elements of ${\cal A}$ are given by 
\begin{eqnarray}
\Lambda^{n}_{~k}(f^{k}_{~m} \star a) &=&(a \star f^{n}_{~k}) \Lambda^{k}_{~m},\\\Lambda^{n}_{~k}(\tilde{f}^{k}_{~m} \tilde{f}^{n}_{~k} \circ S) \Lambda^{k}_{~m} \end{eqnarray}
where now the functionals $f^{n}_{~m}$ and $\tilde{f}^{n}_{~m}$ satisfy 
\begin{eqnarray}
f^{n}_{~m}(ab)=f^{n}_{~k}(a)f^{k}_{~m}(b)  &and~~\tilde{f}^{n}_{~m}(ab)=\tilde{f}^{k}_{~m}(a)\tilde{f}^{n}_{~k}(b)
\end{eqnarray}
for any $a$ and $b \in {\cal A}$. The commutation rules between the elements of ${\cal A}$ and the translations are given by
\begin{eqnarray}
p^{n}a&=&(a \star \tilde{f}^{n}_{~k} \circ S)p^{k} + a \star \tilde{\eta}^{n} \circ S - \Lambda^{n}_{~k}(\tilde{\eta}^{k} \circ S \star a),
\end{eqnarray}
and those between the translations are given by
\begin{eqnarray}
p^{n}p^{m} &=& f^{n}_{~k}(\Lambda^{m}_{~\ell})p^{\ell}p^{k} - (f^{n}_{~k}
(\Lambda^{m}_{~\ell})-\delta^{n}_{~\ell}\delta^{m}_{~k})\tilde{\eta}^{\ell}
(S( \Lambda^{k}_{~p}))p^{p} \nonumber\\&+&\tilde{\eta}^{n}(S(p^{m})) - 
\Lambda^{n}_{~k} \Lambda^{m}_{\ell} \tilde{\eta}^{k}(S(p^{\ell})).
\end{eqnarray}
In addition the generators of ${\cal A}$ sati
sfy also the following relations 
\begin{eqnarray}
S( \Lambda^{k}_{~m}) \Lambda^{a}_{~c} \tilde{\eta}_{k}(p^{c})&=& \tilde{\eta}_{m}(p^{a}),\\\Lambda^{k}_{~m} \tilde{\eta}_{k}(S(p^{a}))&=&\tilde{\eta}_{m}(S(p^{b})) \Lambda^{a}_{~b}.
\end{eqnarray}
Finally it is shown in Ref.[7] that this formalism is consistent if the different functionals  satisfy the following conditions:
\begin{eqnarray}
\tilde{R} = R + I \otimes Q &=& R^{-1} + R^{-1}(I \otimes Q),\\ Q=\lambda
I~~~\lambda \in {\cal C}  &,&~~~\lambda \not= -1,\\
T=-\tilde{R}T~~~ &for&~~\lambda \not = 0\\or~(\tilde{R} + I \otimes I) - 
(\Lambda \Lambda)(\tilde{R}  + I\otimes I)~&for& \lambda =0,\\
(R \otimes I)(I \otimes R)(R \otimes I)&=&(I \otimes R)(R \otimes I)
(I \otimes R),\\
(Z \otimes I)R &+& (\tilde{R} \otimes I)(I \otimes Z)R\nonumber\\ 
=(I \otimes R)(Z \otimes I) &+&(I \otimes R)(\tilde{R} \otimes I)(I \otimes Z),
\\
(R \otimes I - I \otimes I \otimes I)&=&((I \otimes Z)Z -(Z \otimes Z)Z) 
\nonumber\\
+T \otimes I &-& (I \otimes \tilde{R})(\tilde0)
\end{eqnarray}
and
\begin{eqnarray}
(I \otimes R - I \otimes I \otimes I)((\tilde{Z} \otimes I)T - (I \otimes \tilde{Z})T) -(Z \otimes I) - (\tilde{R} \otimes I)(I \otimes Z)T = 0\nonumber\\
\end{eqnarray}
with
\begin{eqnarray}
\tilde{Z}= -RZ
\end{eqnarray}
where $R^{n\ell}_{km}=f^{n}_{~m}(\Lambda^{\ell}_{~k})$, 
$\tilde{R}^{n\ell}_{km} = \tilde{f}^{n}_{~m}(S(\Lambda^{\ell}_{~k}))$, 
$Q^{n}_{k} = \tilde{\eta}_{k}(S(p^{n}))$,
$Z^{nk}_{q} =$ $\tilde{\eta}^{n}(S( \Lambda^{k}_{~q}))$ $= f^{k}_{~q}(S(p^{n}))$, 
$\tilde{Z}^{nk}_{m} = \tilde{f}^{n}_{~m}(S(p^{k}))$ and 
$T^{nm} = \tilde{\eta}^{n}(S(p^{m}))$. 

\section{\bf Lie Algebra of Inhomogeneous Hopf Algebra}
In the following, we shall construct the vector space ${\cal T}$ dual to 
the right-invariant differential one-forms and the quantum commutators of 
different elements of ${\cal T}$.\\ 
We start by\\ 
{\bf Proposition (3,1)}: The one-form\\ 
\begin{eqnarray}
X = \tilde{f}^{a}_{~m}(\Lambda^{m}_{~b})\Theta^{b}_{~a}
\end{eqnarray}
is left- and right-invariant.\\
$Proof.$  X is right-invariant. to show the left-invariance, we replace $a \in {\cal A}$ by $S^{-1}(\Lambda^{a}_{~b})$ into (19) then we apply $S$ on both sides to have 
\begin{eqnarray*}
\Lambda^{c}_{~b}S(\Lambda^{n}_{~k})\tilde{f}^{k}_{~m}(\Lambda^{a}_{~c})= \tilde{f}^{n}_{~k}(\Lambda^{c}_{~b})S(\Lambda^{k}_{~m})\Lambda^{a}_{~c}.
\end{eqnarray*}
Setting the indices $a$ = $m$, we get
\begin{eqnarray}
\Lambda^{c}_{~b}S(\Lambda^{n}_{~k})\tilde{f}^{k}_{~m}(\Lambda^{m}_{~c}) = \tilde{f}^{n}_{~m}(\Lambda^{m}_{~b}).
\end{eqnarray}
Using this relation and (4), we get
\begin{eqnarray}
\Delta_{L}(X) = \Lambda^{c}_{~b}S(\Lambda^{n}_{~k})\tilde{f}^{k}_{~m}(\Lambda^{m}_{~c}) \otimes \theta^{b}_{~n} = \tilde{f}^{n}_{~m}(\Lambda^{m}_{~b}) \otimes \Theta^{b}_{~n} = I \otimes X
\end{eqnarray}
which shows the left-invariance of X.\,\,\,\,\,\,\,\,\,\,\,\, Q.E.D.\\$Remark (3,1)$: Instead of (19), we can use (18) to have an another left- and right-invariant $X=f^{a}_{~m}(S^{-1}(\Lambda^{m}_{~b}))\Theta$: (6-9) show that the classical limit of the functionals are, for any $a \in {\cal B}$, $f^{k}_{~m}(a) = \delta^{k}_{m}\varepsilon(a)$ and $\tilde{f}^{n}_{~k}(a) = \delta^{n}_{k}\varepsilon(a)$ implying $f^{n}_{`m}(S^{-1}(\Lambda^{m}_{~b})) = \delta^{n}_{b}$ and $\tilde{f}^{k}_{~m}(\Lambda^{m}_{~c}) = \delta^{k}_{c}$ and, therefore, (34) can be considered as the quantum trace which reduces to the usual trace in the classical limit. In fact (36) shows the invariance of the quantum trace under the left adjoint coaction action [11].\\We can now express the action of the derivative $d$ on an element of $a \in {\cal B}$ as $da =Xa-aX$ by following the method of differential calculus on the homogeneous quantum groups of Ref.[9]. From the left- and the right-invariance of $X$, it is easy to check the left- and the right-covariance of $d$. In the following we set $V^{a}_{~b} = \tilde{f}^{a}_{~k}(\Lambda^{k}_{~b})$. Using the form of $X$ and (9), we get
\begin{eqnarray}
da &=& Xa-aX\nonumber\\ 
&=& V^{a}_{b}(\Theta^{b}_{a}(a*\varepsilon) - \Theta^{k}_{~\ell}(a*f^{b\ell}_{ak}\circ S) + \Pi^{k} (a*f^{b}_{ak}))
\end{eqnarray}
which can be rewritten under the form
\begin{eqnarray}
da = \Theta^{b}_{~a}(a*L^{a}_{~b}) + \Pi^{a}(a*L_{a})
\end{eqnarray}
where the linear functionnals on ${\cal B}$
\begin{eqnarray}L^{a}_{~b} = V^{a}_{~b}\varepsilon - V^{c}_{~d}(\tilde{f}^{a}_{~c}*f^{d}_{~b})\circ S
\end{eqnarray}and
\begin{eqnarray}
L_{a} = -V^{c}_{~d}(\tilde{\eta}_{c}*f^{d}_{~a}) \circ S 
\end{eqnarray}
are the generators of the Lie algebra of the inhomogeneous Hopf algebra${\cal B}$. Let us note that by virtue of (17)
\begin{eqnarray}
L_{a}(a) = 0  &,~a \in {\cal A},
\end{eqnarray}
and from (37-38), we deduce $L^{a}_{~b}(I) = 0 $ and $L_{a}(I) = 0$.\\
To construct the exterior product of one-forms $\in \Gamma$, we consider, in addition of the right-invariant basis (1), the left-invariant basis $\tilde{\Theta}^{a}_{~b} = S(\Lambda^{a}_{~c}S(\Lambda^{d}_{~b}))\Theta^{c}_{~d}$ and $\tilde{\Pi}^{a} 
+ S(\Lambda^{a}_{~b}S(p^{c}))\Theta^{b}_{~c}$ which transforms as
\begin{eqnarray}
\Delta_{L}(\tilde{\Theta}^{a}_{~b}) &=& I \otimes \tilde{\Theta}^{a}_{~b}~~and~~\Delta_{L}(\tilde{\Pi}^{a}) = I \otimes \tilde{\Pi}^{a},\\\Delta_{R}(\tilde{\Theta}^{a}_{~b}) &=& \tilde{\Theta}^{c}_{~d} \otimes S(\Lambda^{a}_{~c}S(\Lambda^{d}_{~b})),\\\Delta_{R}(\tilde{\Pi}^{a}) &=& \tilde{\Pi}^{b} \otimes S(\Lambda^{a}_{~b}) + \tilde{\Theta}^{b}_{~c} \otimes S(\Lambda^{a}_{~b}S(p^{c})).
\end{eqnarray}
Now we consider the bimodule automorphism $\sigma :\Gamma^{\otimes 2} \rightarrow \Gamma^{\otimes 2}$ such that $\sigma(\tilde{\omega}\otimes _{{\cal B}} \omega) = \omega \otimes _{{\cal B}} \tilde{\omega}$ for any left-invariant element$\tilde{\omega} \in \Gamma$ and any right-invariant element $\omega \in \Gamma$ which gives for the one-forms $\Theta^{a}_{~b}$ 
\begin{eqnarray*}
\sigma(\tilde{\Theta}^{e}_{~f} \otimes _{{\cal B}} \Theta^{g}_{~h}) =S(\Lambda^{e}_{~a}S(\Lambda^{b}_{~f}))\sigma(\Theta^{a}_{~b} \otimes _{{\cal B}} \Theta^{g}_{~h}) = \Theta^{g}_{~h} \otimes _{{\cal B}} S(\Lambda^{e}_{~a}S(\Lambda^{b}_{~f})) \Theta^{a}_{~b}.
\end{eqnarray*}
Using (7) and (17), we get
\begin{eqnarray}
\sigma(\Theta^{a}_{~b} \otimes _{{\cal B}} \Theta^{c}_{~d}) =(\tilde{f}^{g}_{~d}*f^{c}_{~h})(S(\Lambda^{a}_{~e}S(\Lambda^{f}_{~b})))(\Theta^{h}_{~g} \otimes _{{\cal B}} \Theta^{e}_{~f}).
\end{eqnarray}
Similar considerations give 
\begin{eqnarray}
\sigma(\Theta^{a}_{~b} \otimes _{{\cal B}} \Pi^{c}) = f^{c}_{~f} (S(\Lambda^{a}_{~d}S(\Lambda^{e}_{~b})))(\Pi^{f} \otimes _{{\cal B}} \Theta^{d}_{~e})\nonumber\\
 + (\tilde{\eta}^{g} * f^{c}_{~f})(S(\Lambda^{a}_{~d}S(\Lambda^{e}_{~b})))(\Theta^{f}_{~g} \otimes _{{\cal B}} \Theta^{d}_{~e}),\nonumber\\
\sigma(\Pi^{a} \otimes _{{\cal B}} \theta^{b}_{~c}) = (\tilde{f}^{e}_{~c}*f^{b}_{~d})(S(\Lambda^{a}_{~f}))(\Theta^{d}_{~e} \otimes _{{\cal B}} \Pi^{f})\nonumber\\
+(\tilde{\eta}_{c}*f^{b}_{~d})(S(\Lambda^{a}_{~f}S(p^{g})))(\Pi^{d} \otimes _{{\cal B}}^{e}{}_{~c}*f^{b}_{~d})(S(\Lambda^{a}_{~f}S(p^{g})))(\Theta^{d}_{~e} \otimes _{{\cal B}} \Theta^{f}_{~g})\nonumber\\
\end{eqnarray}
and
\begin{eqnarray}
\sigma(\Pi^{a} \otimes _{{\cal B}} \Pi^{b}) &=&f^{b}_{~c}(S(\Lambda^{a}_{~d}))(\Pi^{c} \otimes _{{\cal B}} \Pi^{d})\nonumber\\ + f^{b}_{~c}(S(\Lambda^{a}_{~d}S(p^{e})))(\Pi^{c} \otimes _{{\cal B}} \Theta^{d}_{~e}) &+& (\tilde{\eta}^{d}*f^{b}_{~c})(S(\Lambda^{a}_{~e}))(\Theta^{c}_{~d} \otimes _{{\cal B}} \Pi^{e})\nonumber\\&+& (\tilde{\eta}^{d}*f^{b}_{~c})(S(\Lambda^{a}_{~e}S(p^{f})))(\Theta^{c}_{~d}\otimes _{{\cal B}} \Theta^{e}_{~f}).\nonumber\\
\end{eqnarray} 
In the following we shall write (38) and (45-48) under the compact form
\begin{eqnarray}
da &=& \Pi^{A}(a*L_{A})\\ \sigma(\Pi^{A} \otimes _{{\cal B}} \Pi^{B}) &=& \sigma^{AB}_{~~CD}(\Pi^{C} \otimes _{{\cal B}} \Pi^{D}) 
\end{eqnarray}
where $\Pi^{A} = \Theta ^{a}_{~b}$ and $L_{A} = L^{b}_{~a}$ when the indices $A$ take the values $^{a}_{~b}, ^{c}_{~d},...$ and $\Pi^
{A}= \Pi^{a}$ and $L_{A} = L_{a}$ when the indices $A$ take the values $a, b,...$ and $\sigma^{AB}_{~~CD}$ are given by the coefficients multiplying the tensor products of the one-forms $\Theta$ and $\Pi$ in (45-48). Let us note that from (49), we have $d(a_{(1)})S(a_{(2)} = \Pi^{A}L_{A}(a) = \Pi(a)$ which tranforms as\begin{eqnarray}\Delta_{L}(d(a_{(1)})S(a_{(2)}) = \Delta_{L}(\Pi(a)) =a_{(1)}S(a_{(4)}) \otimes d(a_{(2)})S(a_{(3})\nonumber\\ =a_{(1)}S(a_{(3)} \otimes \Pi(a_{(2)}) =(id \otimes \Pi)Ad(a)\end{eqnarray}where the left adjoint coaction action $Ad:{\cal B} \rightarrow {\cal B} \otimes {\cal B}$ is defined as $Ad(a) = a_{(1)}S(a_{(3)}) \otimes a_{(2)}$ and satisfies 
\begin{eqnarray}
(\Delta \otimes id)Ad = (id \otimes Ad)Ad.
\end{eqnarray}
Here $a_{(1)}$, $a_{(2)}$ and $a_{(3)}$ are element of ${\cal B}$ such that $a_{(1)} \otimes a_{(2)} \otimes a_{(3)} = (id \otimes \Delta)\Delta(a)$.\\  {\bf Theorem (3,1)}: 1. For any $a,b \in {\cal B}$\\
\begin{eqnarray}
L^{a}_{~b}(ab) = L^{a}_{~b}(aa)_{~c}*f^{d}_{b})(S(a))L^{c}_{~d}(b) + (\tilde{\eta}^{a}*f^{c}_{~b})(S(a))L_{c}(b)\\L_{a}(ab)=L_{a}(a)\varepsilon(b) + f^{b}_{~a}(S(a))L_{b}(b) + (\tilde{\eta}_{c}*f^{d}_{~a})(S(a))L^{c}_{~d}(b)
\end{eqnarray}
2. for any $a \in {\cal B}$, the generators $L^{a}_{~b}$ and $L_{a}$ close on the following quantum Lie algebra
\begin{eqnarray}
\left[L^{b}_{~a},L^{d}_{~c}\right](a) &=&(V^{b}_{~a}\delta^{g}_{c}\delta^{d}_{h} - V^{k}_{~n}(\tilde{f}^{b}_{~k}*f^{n}_{~a})(S(\Lambda^{g}_{~c}S(\Lambda^{d}_{~h}))))L^{h}_{~g}(a))\nonumber\\ &-&V^{k}_{~n}(\tilde{f}^{b}_{~k}*f^{n}_{~a})(S(\Lambda^{g}_{~c}S(p^{d})))L_{g}(a),\\\left[L^{b}_{~a},L_{c}\right](a)&=&(V^{b}_{~a}\delta^{g}_{c} -V^{k}_{~n}(\tilde{f}^{b}_{~k}*f^{n}_{~a})(S(\Lambda^{g}_{~c}))L_{g}(a),\\\left[L_{a},L^{c}_{~d}\right](a)&=&V^{k}_{~n}(\tilde{\eta}_{k}*f^{n}_{~a})(S(\Lambda^{g}_{~c}S(p^{d})))L_{g}(a)
\end{eqnarray}
and
\begin{eqnarray}
\left[L_{a},L_{b}\right](a)=0
\end{eqnarray}
satisfying the Jacobi identity
\begin{eqnarray}
\left[L_{A},\left[L_{B},L_{C}\right]\right](a) = \left[\left[L_{A},L_{B}\right],L_{C}\right](a) + \sigma^{EF}_{AB}\left[L_{E},\left[L_{F},L_{C}\right]\right](a)
\end{eqnarray}
$proof.$ Using the Leibniz rule of the derivative $d$, (38) and (8-9), one obtains
\begin{eqnarray*}
d(ab)= \Theta^{b}_{~a}(ab*L^{a}_{~b})&+&\Pi^{a}(ab*L_{a})\\=(\Theta^{~b}_{a}(a*L^{a}_{~b})+\Pi^{a}(a*L_{a}))b &+& a(\Theta^{d}_{~c}(b*L^{c}_{~d})+\Pi^{c}(b*L_{c}))\\=\Theta^{b}_{~a}((a*L^{a}_{~b})b + (a*(\tilde{f}^{a}_{~c}*f^{d}_{~b})\circ S)(b*L^{c}_{~d})&+&(a*\tilde{\eta}^{a}*f^{c}_{~b})\circ S)(b*L_{c}))\\ +\Pi^{a}((a*L_{a})b + (a*f^{b}_{~a} \circ S)(b*L_{b}) &+& (a * (\tilde{\eta}_c*f^{d}_{~a})\circ S)(b*L^{c}_{~d}))
\end{eqnarray*}
for any $a,b \in {\cal B}$. Identifying the coefficients multiplying $\Theta^{b}_{~a}$ and $\Pi_{a}$, then acting $\varepsilon$, one obtains the twisted Leibniz rules for $L^{a}_{~b}$ and $L_{a}$ as (53-54).\\Applying the exterior derivative $d$ on both sides of (49), using $d^{
   2}=0$ 
ade of $d$, we get
\begin{eqnarray}
d^{2}a = 0 =d(\Pi^{A})(a*L_{A}) - (\Pi^{A} \wedge \Pi^{B})(a*(L_{A}*L_{B}),
\end{eqnarray}
for any $a \in{\cal B}$, leading to the Maurer-Cartan formula
\begin{eqnarray}
d\Pi^{A}L_{A}(a) = d\Pi(a) &=& \Pi^{A}  \wedge \Pi^{B}(L_{A}*L_{B})(a)\\&=& (\Pi^{A}\otimes \Pi^{B})[L_{A},L_{B}](a)
\end{eqnarray}
where the exterior product od two one-forms $\Pi^{A}$ and $\Pi^{B}$ is given by $\Pi^{A} \wedge \Pi^{B} = (\Pi^{C}\otimes \Pi^{D})(\delta^{A}_{C}\delta^{B}_{C} - \sigma^{AB}_{~~CD})$ and the commutator $[L_{A},L_{B}](a) = (\delta^{C}_{A}\delta^{D}_{B} - \sigma^{CD}_{~~AB})(L_{C}*L_{B})(a)$ for any $a \in {\cal B}$.\\In the other hand the extension of the derivation (37) to the forms gives
\begin{eqnarray*}
d\Pi^{A}L_{A}(a) = X \wedge \Pi^{A}L_{A}(a) + \Pi^{A}L_{A}(a) \wedge X\\=(I - \sigma)(X \otimes \Pi(a) + \Pi \otimes X).
\end{eqnarray*}
From the left-invariance of $X$ and the property of $\sigma$, it follows
\begin{eqnarray}
d\Pi(a) &=& X\otimes \Pi(a) - \sigma(\Pi(a) \otimes X)\nonumber\\
= (\Theta^{a}_{~b} \otimes \Theta^{c}_{~d} - \sigma (\Theta^{c}_{~d} \otimes \Theta^{a}_{~b}))V^{b}_{~a}L^{d}_{~c}(a) &+& (\Theta^{a}_{~b} \otimes  \Pi^{c} + \sigma(\Pi^{c} \otimes \Theta^{a}_{~b}))V^{b}_{~a}L_{c}(a).\nonumber\\
\end{eqnarray}
We can now use the different stage of the theorem (5,3) of Ref.[9] to get
\begin{eqnarray}
d\Pi^{A}L_{A}(a) = (\Pi^{A} \otimes \Pi^{B})[L_{A},L_{B}](a) = (\Pi^{A} \otimes \Pi^{B})(L_{A} \otimes L_{B})Ad(a)
\end{eqnarray}
Comparing (61) with (63) we get 
\begin{eqnarray*}
\left[L_{A},L_{B}\right](a) = (L_{A} \star L_{B} - \sigma^{CD}_{AB}(L_{C} \star L_{D}))(a) = (L_{A} \otimes L_{B})Ad(a)
\end{eqnarray*}
for any $a \in {\cal B}$. Using (45-48), (63) can be rewritten as
\begin{eqnarray}
d\Pi^{A}L_{A}(a) &=& (\Theta^{a}_{~b} \otimes \Theta^{c}_{~d})((V^{b}_{~a}\delta^{g}_{c}\delta^{d}_{h} - V^{k}_{~n}(\tilde{f}^{b}_{~k}*f^{n}_{~a})(S(\Lambda^{g}_{~c}S(\Lambda^{d}_{~h})))L^{h}_{~g}(a))\nonumber\\ 
&-&V^{k}_{~n}(\tilde\Lambda^{g}_{~c}S(p^{d})))L_{g}(a))\nonumber\\ 
&+& (\Theta^{a}_{~b} \otimes \Pi^{c})(V^{b}_{~a}\delta^{g}_{c} -V^{k}_{~n}(\tilde{f}^{b}_{~k}*f^{n}_{~a})(S(\Lambda^{g}_{~c}))L_{g}(a))\nonumber\\
&-&(\Pi^{a} \otimes \Theta^{c}_{~d})V^{k}_{~n}(\tilde{\eta}_{k}*f^{n}_{~a})(S(\Lambda^{g}_{~c}S(p^{d})))L_{g}(a)).
\end{eqnarray}
Developping now (64) and comparing the coefficients multiplying the  different tensor products of $\Theta^{a}_{~b}$ and $\Pi^{a}$ with the latter relation, we get (55-58).\\Finally we can establish the quantum Jacobi identity (59) by using (52) and (64) as
\begin{eqnarray*}
\left[L_{A},\left[L_{B},L_{C}\right]\right](a) &=& (L_{A} \otimes L_{B} \otimes L_{C})(id \otimes Ad)Ad(a) \\
&=& (L_{A} \otimes L_{B} \otimes L_{C})(\Delta \otimes id)Ad(a) \\
=\left[(L_{A} \star L_{B}),L_{C}\right](a) &=& \left[\left[L_{A},L_{B}\right],L_{C}\right](a) + \sigma^{DE}_{AB}\left[(L_{D} \star L_{E}),L_{C}\right](a)\\
&=& \left[\left[L_{A},L_{B}\right],L_{C}\right](a) + \sigma^{DE}_{AB}\left[L_{D},\left[L_{E},L_{C}\right]\right](a).  
\end{eqnarray*}
\,\,\,\,\,\,\,Q.E.D.\\
As a consequence of the relation (53-54) of this theorem, the universal enveloping algebra ${{\cal U}(\cal T})$ can be equipped with a hopf algebra structure whose the comultiplication is given by
\begin{eqnarray*}
\Delta'(L^{a}_{~b}) = L^{a}_{~b} \otimes \varepsilon + S'(f^{da}_{cb})\otimes L^{c}_{~d} + S'(f^{ca}_{~b})\otimes L_{c},\\
\Delta'(L_{a}) = L_{a} \otimes \varepsilon + S'(f^{b}_{~a}) \otimes L_{b} + S'(f^{d}_{ca}) \otimes L^{c}_{d},
\end{eqnarray*}
the andipode $S'$ by
\begin{eqnarray*}
S'(L^{a}_{~b}) = -S'^{2}(f^{da}_{cb})\star L^{c}_{~d} - S'^{2}(f^{ca}_{~b})\star L_{c},\\S'(L_{a}) = - S'^{2}(f^{b}_{~a})\star L_{b} - S'^{2}(f^{d}_{ca})\star L^{c}_{~d}
\end{eqnarray*}
 and the counit by
\begin{eqnarray*}
\varepsilon'(L^{a}_{~b}) = \varepsilon'(L_{b}) =0.
\end{eqnarray*}
Using the Hopf algebra structure of ${\cal B}^{0}$, one can see that the axioms of the Hopf algebra are satisfied for ${\cal U}$(in
 virtue of (41)), the relations (55-58) reduce to the homogeneous part of (55) which is the Lie bracket verified by the generators of  the Lie  algebra of the homogeneous quantum groups [10-11]. For any $a \in {\cal B}$ of the form ${\cal A}p^{n} + {\cal A}$ the convolution product $(L_{a} \star L_{b})(a)=0$, due to (41), hence (58) is non-trivial only when it apply on an element of ${\cal B}$ which contains products of two translations and more.\\

{\bf Acknowledgments:} We are grateful to M. Dubois-Violette for valuable discussions. The author (M. L.) wishes to thank the members of the L.P.T.H.E. of Orsay University for their hospitality during November 1997 where a part of the work of this paper was done.\\

{\bf References:}\\
1-- O. Ogiesvestsky, W. B. Schmidke, J. Wess and B. Zumino, Commun. Math. Phys. 150(1992)495.\\
2-- M. Schlicker, W. Weich and R. Weixler, Lett. Math. Phys. 27(1993)217.\\
3-- V. K. Dobrev, J. Phys. A:Math. Gen. 26(1993)1317.\\
4-- J. Lukierski, A.  Nowicki, H. Ruegg and V. N. Tolstoy, Phys. Lett. B264
(1991)331.\\
5-- J. Lukierski and A. Nowicki, Phys. Lett. B279(1992)299.\\
6-- L. Castellani, Lett. Math. Phys. 30(1994)233, Commun. Math. Phys. 171(1995)
383, P. Aschieri and L. Castellani, Int. J. Mod. Phys. A11(1996)4513\\
7-- M. Lagraa and N. Touhami: "The noncommutative Hopf algebra", Preprint L.P.T.O. Es-Senia, q-alg 9705005 and L.P.T.H.E.-Orsay 97/51 available from http://qcd.th.u-psud.fr .\\
8-- P. Podles and S. L. Woronowicz, Commun. Math, Phys. 185(1997)325.\\
9-- S. L. Woronowicz, Commun. Math. Phys. 122(1989)125.\\
10-- U. Carow-Watamura, M. Schliecker, S. Watamura and W. Weich, Commun. Math. 
Phys. 142(1991)605.\\
11-- M. Lagraa, Int. J. Mod. Phys. A11(1996)699. \end{document}